\documentclass[12pt]{amsart}
\usepackage{amssymb, amscd, amsmath, amsthm, 
epsf, epsfig, latexsym, enumerate}

\begin{document}
 
\title{A small infinitely-ended 2-knot group}

\author{R. Budney \and J. A. Hillman }
\address{{Mathematics and Statistics, Victoria University,}
\newline
{Victoria, V8W 3R4, Canada} 
\newline
{School of Mathematics and Statistics, University of Sydney,}
\newline
Sydney,  NSW 2006, Australia }

\email{rybu@uvic.ca,jonathan.hillman@sydney.edu.au}
 
\begin{abstract} 
We show that a 2-knot group discovered in the course of a census of 4-manifolds 
with small triangulations is an HNN extension with finite base and proper associated subgroups,
and has the smallest base among such knot groups.
\end{abstract}


 
\maketitle 

\section{introduction}

Nontrivial classical knot groups have one end. 
This is equivalent to the asphericity of the knot complement, 
in the light of Dehn's Lemma and Poincar\'e duality in the universal cover.
In higher dimensions the complements of nontrivial knots are never aspherical.
However, Kervaire gave a uniform algebraic characterization of $n$-knot groups
for all $n\geq3$, and a partial characterization of 2-knot groups \cite{Ke}.
(Artin's spinning construction shows that classical knot groups are 2-knot groups
and 2-knot groups are high dimensional knot groups.)
These characterizations have been used to provide examples of such groups 
with various  properties.
In particular, they have been used to find knots whose groups have more than one end.

The 2-twist spins $\tau_2K$ of 2-bridge knots $K$ provide many examples 
of 2-knot groups with two ends. 
The first examples of higher dimensional knots
with infinitely-ended groups were given in \cite{GAM78}.
Their examples are 2-knots, 
and the groups are HNN extensions with finite base and proper associated subgroups.
The simplest of these has presentation
\[
\langle{a,b,t}\mid{tat^{-1}=a^2},~a^3=1,~aba^{-1}=b^2\rangle,
\]
with base $\langle{a,b}\rangle\cong{Z/7Z}\rtimes_2Z/3Z$ and 
both associated subgroups $\langle{a}\rangle\cong{Z/3Z}$.
(Note that the second and third relations together imply that $b^7=1$.)
It is clear from this presentation that the group is also a free product with amalgamation 
of $\langle{a,b}\rangle$ with $\langle{a,t}\rangle\cong\pi\tau_23_1$ 
over $\langle{a}\rangle$, where $3_1$ is the trefoil knot \cite{Ra81}.
In fact it is the group of a satellite of $\tau_23_1$ with companion Fox's Example 10,
as is clear from the analysis of the groups of such knots in \cite{Ka83}.
Replacing the trefoil  with other 2-bridge knots gives all of the examples of \S1 of \cite{GAM78}.

In \S2 we shall show that a knot exterior recently discovered in the course 
of a computational census of 4-manifolds with small triangulations has group $\pi$ 
which is an HNN extension with base the generalized quaternionic group
$Q(16)$ and (distinct) associated subgroups $Q(8)$.
This group  is not properly the group of a satellite knot (\S3).
In \S4 we show that no smaller finite group is the base of an HNN extension which
is an infinitely-ended knot group, and that there is only one such group with base $Q(16)$.
However  there are such groups with base $D_8\times{Z/2Z}$
and associated subgroups $(Z/2Z)^3$ (\S5).
In \S6 we show that all knots with the exterior described in \S2 are strongly $+$amphicheiral,
and we determine generators for $Out(\pi)\cong(Z/2Z)^3$.

For us, an $n$-knot is a locally flat embedding of $S^n$ in $S^{n+2}$, 
and so orientations of the spheres determine orientations 
of the knot complements and preferred meridians for the knots.
We shall use Chapter 14 of \cite{Hi} as a one-stop reference for the aspects of
higher dimensional knot theory that we need.
All homology groups considered below shall have integral coefficients,
and so we shall write $H_i(G)$ instead of $H_i(G;\mathbb{Z})$.

\section{The knot exterior}\label{MANDETAILS}

The knot exterior of this paper was discovered while forming a census of $4$-manifolds
triangulable with $6$ or less pentachora.  
Precisely, Ben Burton generalized the census-generation algorithm in Regina \cite{REG} 
to enumerate all unordered $4$-dimensional delta-complexes whose vertex links 
are triangulated $3$-spheres or more generally $3$-manifolds. 
Triangulations having a non-spherical manifold vertex link are called {\it cusped triangulated manifolds}. Most of the non-trivial knot exteriors in the census are of ideal/cusped type.  
A previous paper was written on the simplest non-trivial knot exterior in the census \cite{BBH}.  
A future paper will describe the census in full.  

In the census there are approximately 1.4 million combinatorial classes of knot exteriors 
in homotopy spheres.  By {\it combinatorial class} we mean triangulated manifolds, 
up to homeomorphism that preserve the simplicial subdivision from the triangulation, 
i.e. the homeomorphisms need not preserve the characteristic maps of the individual simplices. 
Among these 1.4 million triangulations, 8521 have non-abelian fundamental group.
Most of these have finitely generated commutator subgroup.
There are just twenty exceptional cases.

Ten of these have group $\Phi$ with presentation $\langle{a,t}\mid tat^{-1}=a^2\rangle$, 
and so are homeomorphic to the exterior of Fox's Example 10.
The final ten all have fundamental group isomorphic to the group $\pi$
with presentation
\[
\langle{a,t}\mid{a^8},
~a^{-2}tata^2t^{-2},~a^2tata^{-2}t^{-2},~a^4t^{-1}a^{-4}t\rangle.
\]

We briefly describe one of the ten triangulations giving a manifold $M$ with fundamental group 
$\pi$. At present we know there are at most two PL homeomorphism types represented by these triangulations.
\vskip 0.5cm
\centerline{
\begin{tabular}{l | c c c c c}
M &  (0123) &  (0124) &  (0134) &  (0234) &  (1234) \\
\hline
0 & 4 (0123) & 3 (0124) & 2 (1320) & 2 (0234) & 1 (1234) \\
1 & 5 (0123) & 4 (1240) & 4 (4320) & 2 (1234) & 0 (1234) \\
2 & 0 (4031) & 5 (3142) & 3 (4013) & 0 (0234) & 1 (0234) \\
3 & 5 (0421) & 0 (0124) & 2 (1340) & 4 (1423) & 5 (4031) \\
4 & 0 (0123) & 1 (4012) & 5 (0324) & 1 (4310) & 3 (0342) \\
5 & 1 (0123) & 3 (0321) & 3 (2431) & 4 (0314) & 2 (1402) \\
\hline
\end{tabular}
}
\vskip 0.5cm
The leftmost column lists the pentachora of the triangulation, labelling them $0$ through $5$. 
In each row, to the right of the pentachoron index is a collection of pairs $n (abcd)$.  The first
number $n$ indicates on which pentachoron this tetrahedral facet is glued to. 
The entry $(abcd)$ indicates the affine-linear map on the tetrahedral facet.  
For example, the row $1 \ 5 (0123) \ 4 (2104) \ 4 (4320) \ 2 (1234) \ 0 (1234)$ indicates that
in the pentachoron indexed by $1$, the $4$th tetrahedron is glued to the $4$th tetrahedron but in pentachoron $5$, with vertices $(0123)$ glued to  $(0123)$ in that order.  Similarly, the $2$nd tetrahedron is glued to the $1$st tetrahedron but in pentachoron $4$, with vertices $(0134)$ sent to $(4320)$ in that order, etc.  We leave the readers to consult the documentation for \cite{REG} for details.  In summary, this triangulation has no internal vertices -- after performing the above gluings, all the vertices have been identified, thus the single vertex has vertex link a $30$-tetrahedron triangulation of $S^1 \times S^2$.  The triangulation has $3$ edges, $12$ triangles, $15$ tetrahedra and $6$ pentachora. 

Verification that $M$ is a knot exterior in a homotopy sphere is similar to the argument 
in \cite{BBH} and left to the reader.  
(Up to changes of orientation, there are at most two knots with exterior homeomorphic to $M$.)
Budney and Burton have automated the process and it is implemented in the software \cite{REG}. Perhaps for some readers it would be more appealing to read the algorithm implemented in Regina.  At present the $4$-manifolds tools are in the {\it development repository} of Regina, and these tools will be in the general release of Regina by version 5.0.
One can readily check that the above triangulation has a single non-trivial symmetry,
an involution that reverses orientation and acts non-trivially on $H_1(M)$.
The involution is the simplicial map that sends pentachoron 0 to pentachoron 2,
sending vertices $(01234)$ to $(24103)$, in that order.
Pentachoron 1 is sent to pentachoron 3, 
sending vertices $(01234)$ to $(23041)$, in that order, 
and pentachoron 4 to pentachoron 5, sending vertices $(01234)$ to $(42130)$, in that order.

Returning to the group $\pi=\pi_1(M)$, 
we find that the third relator in the above presentation 
is a consequence of the others, 
and so this presentation simplifies to
\[
\langle{a,t}\mid{tat^{-1}=a^2t^2a^{-2}t^{-2}},~a^8,~a^4t=ta^4
\rangle.
\]

Setting $b=ta^2t^{-1}$, 
this becomes
\[
\langle{a,b,t}\mid{ta^2t^{-1}=b},~tabt^{-1}=a^2,~a^4=b^2=(ab)^2
\rangle.
\]
(The final two relations imply that $bab^{-1}=a^{-1}$. 
Hence $a^8=1$, and so $a^4$ is a central involution.)
Thus $\pi\cong{B*_H\varphi}$ is an HNN extension, 
with base $B\cong{Q(16)}$, 
the generalized quaternionic group with presentation
\[
\langle{a,b}\mid{a^4=b^2=(ab)^2}\rangle.
\]
and associated subgroups 
$H=\langle{a^2,ab}\rangle\cong\varphi(H)=\langle{a^2,b}\rangle\cong{Q(8)}$.
The commutator subgroup is an iterated generalized free product with amalgamation
\[
\pi'\cong\dots B*_HB*_HB\dots,
\]
and is perfect ($\pi'=\pi''$).
Since $H$ and $\varphi(H)$ are proper subgroups of $B$ the commutator subgroup
is not finitely generated.
Hence no knot with this group is fibred.

\section{$\pi$ is not properly the group of a satellite knot}

If an $n$-knot $K$ is a satellite of $K_1$ about $K_2$ relative to a simple closed curve 
$\gamma$ in $X(K_1)$ then 
\[
\pi{K}\cong\pi{K_2}/\langle\langle{w^q}\rangle\rangle*_{w=[\gamma]}\pi{K_1},
\]
where $[\gamma]\in\pi{K_1}$ has order $q\geq0$ and $w$ is a meridian for $K_2$.
The case $q=0$ corresponds to $[\gamma]$ having infinite order.
(See \cite{Ka83}, or page 271 of \cite{Hi}.)
If $K_2=\tau_rk$ is the $r$-twist spin of an $(n-1)$-knot $k$ and $(q,r)=1$ then
$\pi{K_2}/\langle\langle{w^q}\rangle\rangle\cong{Z/qZ}$.
Thus every 2-knot group with non-trivial torsion is trivially the group of a satellite knot.
We shall say that $\pi{K}$ is {\it properly\/} the group of a satellite knot if 
$|\pi{K_2}/\langle\langle{w^q}\rangle\rangle|>q$.

Suppose that the group $\pi$ of \S2 is properly the group of a satellite knot.
Since  $\pi$  has a central subgroup of order 2,  it is not of the form $A*_\mathbb{Z}B$
with $A$ and $B$ nontrivial knot groups.
Hence $\pi\cong{G}*_{Z/qZ}H$,
where $G$ is a knot group and $H$ is the quotient of a knot group by the $q$th power of a meridian, for some $q>0$, but is not cyclic.

Since $\pi$ is an HNN extension with finite base it is virtually free.
(See Corollary IV.1.9 of \cite{DD}.
The argument given there implies that $\pi$ has a free subgroup of index dividing $16!$.) 
Therefore so are $G$ and $H$, 
and all these groups have well-defined virtual Euler characteristics.
Mayer-Vietoris arguments give
\[
\chi^v(\pi)=\chi^v(Q(16))-\chi^v(Q(8))=\frac1{16}-\frac18=-\frac1{16}
\]
and
\[
\chi^v(\pi)=\chi^v(G)+\chi^v(H)-\frac1q,
\]
while $\chi^v(G)\leq0$, since $G$ is an infinite virtually free group.
Hence
\[
\chi^v(H)\geq\frac1q-\frac1{16}.
\] 

Now since $\pi\cong{Q(16)}*_{Q(8)}\varphi$, any finite subgroup of $\pi$ is conjugate 
to a subgroup of $Q(16)$.
Therefore $q$ divides 8, and so $\chi^v(H)>0$.
Hence $H$ is finite, and so it is isomorphic to a subgroup of $Q(16)$.

We then find that the only possibility is that $q=8$ and $H\cong{Q(16)}$.
But then $\chi^v(G)=0$, and so $G'$ is finite, 
and is either $Z/8Z$ or $Q(16)$.
Neither of these groups admits a meridianal automorphism, 
and so there is no such knot group $G$.
Thus we may conclude that $\pi$ is not properly the group of a satellite knot.

HNN extensions arise naturally in knot theory when an $n$-knot $K$ has a minimal 
Seifert hypersurface $V$, one for which the pushoffs from $V$ to either side of 
$Y=S^{n+2}\setminus{V}$ are both $\pi_1$-injective.
The knot group $\pi{K}$ is then an HNN extension with base $\pi_1(Y)$ and associated subgroups isomorphic to $\pi_1(V)$.
There are 2-knots  with group ${Z/3Z}\rtimes\mathbb{Z}$ 
(the group of the 2-twist spun trefoil) which do not have minimal Seifert hypersurfaces.
(See Chapter 17 of \cite{Hi}.)
Does a 2-knot with exterior the manifold $M$ of \S2 have a minimal Seifert 
hypersurface realizing the HNN structure $\pi\cong{Q(16)}*_{Q(8)}\varphi$?
(Knots  related by composition with reflections of $S^n$ or $S^{n+2}$ 
have similar Seifert hypersurfaces.)

\section{hnn extensions with small finite base}

Let $G=B*_C\phi$ be an HNN extension with base $B$
and associated subgroups $C$ and $\phi(C)$.
Let $j:C\to{B}$ be the natural inclusion.
Consideration of the Mayer-Vietoris sequence for the extension shows that 
$H_1(G)\cong\mathbb{Z}$ and $H_2(G)=0$ if and only if 
$H_1(j)-H_1(\phi)$ is an isomorphism and $H_2(j)-H_2(\phi)$ is surjective. 
If the $H_1$ condition holds then $B/N$ is perfect,
where $N=\langle\langle\{j(g)^{-1}\phi(g)|g\in{C}\}\rangle\rangle$ 
is the normal closure of  $\{j(g)^{-1}\phi(g)|g\in{C}\}$ in $B$.
In particular, if $B$ is solvable then $N=B$.
The $H_2$ condition holds automatically if $H_2(B;\mathbb{Z})=0$,
in particular, if $B$ is a finite group of cohomological period 4.

If both homological conditions hold and $N=B$ then the stable letter 
of the HNN extension normally generates $G$,
and so $G$ is a knot group \cite{GAM78}.
(However, such HNN extensions need not be 2-knot groups.
The group $Z/5Z\rtimes_2Z$ is a high dimensional knot group which is
an HNN extension with $H=B=Z/5Z$, but the Farber-Levine condition fails,
since $ 2^2\not\equiv1$ {\it mod} (5).
See Chapter 14 of \cite{Hi}.)

In particular, if $B$ is finite and $C$ is a proper subgroup 
then $B$ is nonabelian, and $C/C'\cong{B/B'}$.
Hence $H_1(B)$ cannot be cyclic of even order.
Therefore $B$ is neither a dihedral group $D_{2k}$ with $k$ odd,
nor $Z/3Z\rtimes_{-1}Z/4Z$. 
This leaves only $Q(8)$, $D_8$, $D_{12}$ and $A_4$
among the groups of order less than 16. 

The group $Q(8)$ has no proper subgroup with abelianization $(Z/2Z)^2$.
If $B=D_{4k}$ then $C$ must be isomorphic to $D_{4l}$, 
for some $l$ dividing $k$.
But then $H_2(C)\cong{H_2(B)}=Z/2Z$,
and $H_2(j)$ and $H_2(\phi)$ are the same isomorphism.
Hence $H_2(j)-H_2(\phi)$ is not an epimorphism.
Thus we may exclude $D_8$ and $D_{12}$.
If $B=A_4$ then $C$ must be $Z/3Z$. 
Since $H_2(Z/3Z)=0$ and $H_2(A_4)=Z/2Z$,
this group may be excluded also.

Thus the smallest possible base must have order at least 16.
The group $Q(16)$ has two proper subgroups with abelianization $(Z/2Z)^2$.
These are $\langle{a^2,b}\rangle$ and $\langle{a^2,ab}\rangle$,
and are isomorphic to $Q(8)$.
The automorphism $a\mapsto{a},b\mapsto{ab}$ of $Q(16)$
carries one onto the other.
Fix generators $x,y$ for $Q(8)$. 
Then we may assume that  $j(x)=a^2$ and $j(y)=b$.
If $\phi$ is another embedding then $\phi(x)$ and $\phi(y)$ have order 4,
so one must be $a^{\pm2}$ and the other $a^ib$.
If, moreover, $H_1(j)-H_1(\phi)$ is an isomorphism then $\phi(x)=a^ib$
and $\phi(y)=a^{\pm2}$, and $i$ must be odd.
After conjugation in $Q(16)$ we may assume that $\phi(x)=ab$ and $\phi(y)=a^2$.
Since $Q(16)$ is solvable and $H_2(Q(16))=0$ the HNN extension 
$\pi=Q(16)*_{Q(8)}\varphi$ is a knot group, 
and it has the smallest finite base among all such HNN extensions.
Moreover, it is the unique such group with HNN base $Q(16)$.

\section{further examples with base of order 16}

There are eight other non-abelian groups of order 16.
Four are semidirect products $K\rtimes{L}$ with $K$ and $L$ cyclic.
In three of these four cases $H_2(C)=H_2(B)=Z/2Z$, 
and so these may be ruled out,
by the argument that excluded $D_8$ and $D_{12}$.
The fourth  group $M_{16}$ has presentation
$\langle{a,x}\mid{a^8=x^2=1~,xax=a^5}\rangle$.
The abelianization is $Z/4Z\oplus{Z/2Z}$, and so $C$ must be $\langle{a^2,x}\rangle$.
There is no second embedding $\phi$ such that $H_1(j)-H_1(\phi)$ is an isomorphism, 
and so we may rule out $M_{16}$.
(Note that $H_2(M_{16})=0$.)

The next group to consider is ${(Z/2Z)^2\rtimes_\theta{Z/4Z}}$,
the semidirect product with action generated by  
$\theta=\left(\begin{smallmatrix}
1&1\\ 
0&1
\end{smallmatrix}\right)\in{GL(2,\mathbb{F}_2)}$,
and with presentation 
$\langle{ a,b,x}\mid{a^4=b^2=x^2=1,~ab=ba, ~bx=xb, ~axa^{-1}=bx}\rangle$.
(Since $b=a^2(ax)^2$,  this group is generated by $\{a,x\}$.)
This has abelianization $Z/4Z\oplus{Z/2Z}$, and so $H_2(C)=Z/2Z$.
It follows from the LHS spectral sequence for $B$ as a semidirect product
that $H_2(B)$ maps onto $H_1(Z/4Z;(Z/2Z)^2)=\mathrm{Ker}(\theta-I)=Z/2Z$.
Hence either $H_2(B)=Z/2Z$ and $H_2(j)-H_2(\phi)=0$,
or $H_2(B)$ has order $\geq4$.
In neither case is $H_2(j)-H_2(\phi)$ an epimorphism,
and so we may exclude this group.

The remaining three have abelianization $(Z/2Z)^3$, so
 $C$ must be an abelian subgroup of index 2, and hence normal.
 These are $Q(8)\times{Z/2Z}$, $D_8\times{Z/2Z}$ and 
 the central product of $D_8$ with $Z/4Z$,
 with presentation
 \[
  \langle{ a,c,x}\mid{ a^4=x^2=1,~a^2=c^2,~ac=ca,~cx=xc,~xax=a^{-1}}\rangle.
  \]
We may eliminate $Q(8)\times{Z/2Z}$ and the central product immediately,
as neither has a proper subgroup with abelianization $(Z/2Z)^3$.

The final group is $B=D_8\times{Z/2Z}$, with presentation
\[
\langle{a,b,x}\mid{a^4=b^2=x^2=1,~ab=ba,~bx=xb,~xax=a^{-1}}\rangle.
\]
There are two proper subgroups isomorphic to $(Z/2Z)^3$.
These are $\langle{a^2,b,x}\rangle$ and $\langle{a^2,b,ax}\rangle$,
and the automorphism $a\mapsto{a}$, $b\mapsto{b}$, $x\mapsto{ax}$ of $B$ 
carries one onto the other
Let $\{c_1,c_2,c_3\}$ be the standard basis for $(Z/2Z)^3$.
Then we may assume that $C$ is the image of the embedding $j:(Z/2Z)^3\to{B}$ 
determined by $j(c_1)=a^2$, $j(c_2)=b$ and $j(c_3)=x$.

Let $V=(Z/2Z)^2$, and let $e_i$ be the image of the generator of $H_2(V)=Z/2Z$
under the inclusion of $V$ onto the subgroup generated by $\{c_k|k\not=i\}$.
Then $\{e_1,e_2,e_3\}$  is a basis for $H_2(C)\cong(Z/2Z)^3$.
We also have $H_2(B)\cong(Z/2Z)^3$,
since $H_2(B)=H_2(D_8)\oplus(H_1(D_8)\otimes{Z/2Z})$,
by the K\"unneth Theorem.
This has a basis $\{f_1,f_2,f_3\}$,
where $f_1$ is the image of the generator of $H_2(D_8)$,
$f_2=a\otimes{b}$ and $f_3=x\otimes{b}$.
The homomorphism $H_2(j)$ sends $e_1, e_2$ and $e_3$ to $f_3,f_1$ and 0, respectively.

Reimbeddings satisfying the $H_1$ condition must carry $C$ to 
$\widetilde{C}=\langle{a^2,b,ax}\rangle$.
There are $|GL(3,\mathbb{F}_2)|=168$ possible isomorphisms $\phi$.
Conjugation in $B$ reduces this by half (since $[B:C]=2$), 
but this still leaves too many possibilities to examine easily by hand.
We shall just give one example.

Let  $\tilde{j}(c_1)=ax$, $\tilde{j}(c_2)=a^2$ and $\tilde{j}(c_3)=b$.
Then $\mathrm{Im}(\tilde{j})=\widetilde{C}$, and $H_2(\tilde{j})$ sends $e_1, e_2$ 
and $e_3$ to $0, f_2+f_3$ and $f_1$, respectively.
It follows easily that the homological conditions are satisfied.
Let $\phi=\tilde{j}j^{-1}$ (so $\phi(a^2)=ax$, $\phi(b)=a^2$ and $\phi(x)=b$),
and let $\pi=B*_C\phi$. 
Then the stable letter of the HNN extension is a normal generator for $\pi$, 
since $B$ is solvable, and so $\pi$ is a high-dimensional knot group.
Is there a 2-knot group with HNN basis $B=D_8\times{Z/2Z}$?

\section{symmetries} 

Four of the ten (ideal) triangulations in the census which realize $\pi$ have 
simplicial involutions which reverse orientation and acts non-trivially on $H_1(M)$.
In particular, the triangulation displayed in \S2 has a single non-trivial symmetry,
which is such an involution. 
It is the simplicial map that sends pentachoron 0 to pentachoron 2,
sending vertices $(01234)$ to $(24103)$, in that order.
Pentachoron 1 is sent to pentachoron 3, 
sending vertices $(01234)$ to $(23041)$, in that order, 
and pentachoron 4 to pentachoron 5, sending vertices $(01234)$ to $(42130)$, in that order.

The open 4-manifold $M$ is the interior of a compact 4-manifold $\overline{M}$,
with boundary $\partial\overline{M}\cong{S^1}\times{S^2}$, and the involution 
extends to $\overline{M}$.
It is known that there are 13 involutions of $S^1\times{S^2}$, 
up to conjugacy \cite{To73}.
In Tollefson's classification precisely three reverse both the orientation and the meridian, 
and they are determined by their fixed-point sets. 
One has fixed-point set $S^0\amalg{S^2}$, and does not extend across $D^2\times{S^2}$.
The others have fixed-point set $S^2\amalg{S^2}$ and $S^0\amalg{S^0}$, respectively,
and extend to involutions of $D^2\times{S^2}$.
Computation shows that the present involution fixes $S^0\amalg{S^0}$
(i.e., four points) in $\partial\overline{M}$, 
and thus extends across any homotopy 4-sphere of the form $M\cup{D^2}\times{S^2}$. 
Hence every knot with exterior $M$ is strongly $+$amphicheiral.
Is any such knot also invertible or reflexive?

On the algebraic side, it is easy to determine the outer automorphism classes of $\pi$,
since every automorphism of an HNN extension with finite base must carry the base 
to a conjugate of itself. 
Thus $Out(\pi)$ is generated by automorphisms which fix $Q(16)$ set-wise.
If $\alpha$ is such an automorphism then $\alpha(a)=a^i$ and  $\alpha(b)=a^jb$, 
for some odd $i=\pm1$ and some $j$, and $\alpha(t)=wt^\epsilon$, 
for some $w\in\pi'=\langle\langle{a}\rangle\rangle$ and $\epsilon=\pm1$.

Suppose first that $\epsilon=1$.
The images of the relations give equations
\[
wb^iw^{-1}=a^jb\quad\mathrm{ and}\quad{wta^{i+j}bt^{-1}w^{-1}=a^{2i}}.
\]
The first equation implies that $w\in{Q(16)}$,
by the uniqueness of normal forms for elements of an HNN extension.
Hence $j$ must be even, and so $i+j=2k+1$, for some $k$.
The second equation then becomes 
\[wb^ka^2w^{-1}=a^{2i},
\]
and so $k$ is also even.
On following this through, we find that there is an unique such automorphism
for each $w\in{Q(16)}$.
Four of these are inner automorphisms, given by conjugation
by elements of the subgroup $\langle{a^2,ab}\rangle$,
and so we need only consider the automorphisms $f$ and $g$,
given by $f(a)=a$, $f(b)=b$ and $f(t)=a^4t$,
and $g(a)=a^{-1}$, $g(b)=a^{-2}b$ and $g(t)=at$.
It is easy to see that $fg=gf$ and $f^2=g^2=id_\pi$.

There is also an automorphism $h$ such that $h(a)=a$, $h(b)=ab$ and $h(t)=(at)^{-1}$.
This automorphism induces the involution of $\pi/\pi'\cong\mathbb{Z}$.
We have  $fh=hf$ and $(gh)^2=id_\pi$,
while $h^2$ is conjugation by $a$, and so $h^8=id_\pi$.
Thus $Out(\pi)\cong(Z/2Z)^3$, and is generated by the images of $f$, $g$ and $h$.

The HNN structure determines a  Mayer-Vietoris sequence
\[
\dots\to{H_4(Q(16))}\to{H_4(\pi)}\to{H_3(Q(8))}\to{H_3(Q(16))}\to\dots,
\]
where the right hand homomorphism is the difference
of the homomorphisms induced
by the inclusions of the two associated subgroups.
Since $H_3(Q(8))\cong{Z/8Z}$, $H_3(Q(16))\cong{Z/16Z}$ and $H_4(Q(16))=0$, it
follows that $H_4(\pi)$ is cyclic of order dividing 8.
Since both of the homomorphisms $H_3(Q(8))\to{H_3(Q(16))}$ induced by the inclusions are 
injective, $H_3(\pi)\not=0$.
The automorphisms $f$ and $g$ preserve the associated subgroups,
and induce the identity on $H_3(Q(8))$.
Hence they also induce the identity on $H_4(\pi)$. 
How does $h$ act on this homology group?



\begin{thebibliography}{99}


\bibitem{BBH} Budney, R., Burton, B. and Hillman, J. A. 
Triangulating a Cappell-Shaneson knot complement,
Math. Res. Lett. 19 (2012), 1117--1126.

\bibitem{REG} Burton, B. A.,  Budney, R., Pettersson, W. et al.,
Regina: Software for 3-manifold topology and normal surface theory,

http://regina.sourceforge.net/, 1999--2013.

\bibitem{BZ} Burde, G. and Zieschang, H. {\it Knots},

de Gruyter Studies in Mathematics 5,

W. de Gruyter Verlag, Berlin -- New York (1985).

\bibitem{DD} Dicks, W. and Dunwoody, M. J. {\it Groups Acting on Graphs},

Cambridge Studies in Advanced Mathematics 17,

Cambridge University Press, Cambridge -- New York -- Melbourne (1989).

\bibitem{GAM78} Gonz\'alez-Acu\~na, F. and Montesinos, J. M.
Ends of knot groups,

Ann. Math. 108 (1978), 91--96.

\bibitem{Hi} Hillman, J. A. {\it Four-Manifolds, Geometries and Knots},

Geometry and Topology Monographs 5, 

Geometry and Topology Publications (2002). (Revision 2007).

\bibitem{Ka83} Kanenobu, T. Groups of higher-dimensional satellite knots,

J. Pure Appl. Alg. 28 (1983), 179--188.

\bibitem{Ke} Kervaire, M. A. Les noeuds de dimensions sup\'erieures,

Bull. Soc. Math. France 93 (1965), 225--271.

\bibitem{Ra81} Ratcliffe, J. On the ends of higher dimensional knot groups,

J. Pure Appl. Alg. 20 (1981), 317--324.

\bibitem{To73} Tollefson, J. L. Involutions on $S^1\times{S^2}$ and other manifolds,

Trans. Amer. Math. Soc. 183 (1973), 139--152.

\end{thebibliography}
\end{document}